\documentclass{amsart}
\usepackage{amssymb}

\newcommand{\jj}{\vee}% join
\newcommand{\mm}{\wedge}% meet
\newcommand{\JJ}{\bigvee}% big join
\newcommand{\JJm}[2]{\JJ(\,#1\mid#2\,)}% big join with a middle
\newcommand{\uu}{\cup}% union
\newcommand{\ii}{\cap}% intersection
\newcommand{\UU}{\bigcup}% big union
\newcommand{\II}{\bigcap}% big intersection
\newcommand{\UUm}[2]{\UU(\,#1\mid#2\,)}% big union with a middle
\newcommand{\IIm}[2]{\II(\,#1\mid#2\,)}% big intersection with a middle

\newcommand{\ci}{\subseteq}% contained in with equality
\newcommand{\nin}{\notin}% not \in
\newcommand{\es}{\varnothing}% the empty set
\newcommand{\set}[1]{\{#1\}}% set 
\newcommand{\setm}[2]{\{\,#1\mid#2\,\}}% set with a middle
\def\vv<#1>{\langle#1\rangle}% vector 

\newcommand{\ga}{\alpha}
\newcommand{\gb}{\beta}
\renewcommand{\ge}{\varepsilon}% use \geq for >=
\newcommand{\gf}{\varphi}
\newcommand{\gm}{\mu}
\newcommand{\go}{\omega}
\newcommand{\gx}{\xi}
\newcommand{\gy}{\psi}

\newcommand{\gQ}{\Theta}

\newcommand{\tbf}{\textbf}% text bold

\newcommand{\E}[1]{\mathcal{#1}}% same as \C
\newcommand{\ol}[1]{\overline{#1}}
\def\con#1=#2(#3){#1\equiv#2\pod{#3}}
\providecommand{\bysame}{\makebox[3em]{\hrulefill}\thinspace}
\newcommand{\q}{\quad}

\newcommand{\iso}{\cong}

\theoremstyle{plain}
\newtheorem{theorem}{Theorem}
\newtheorem{lemma}{Lemma}[section]

\newtheorem{corollary}[lemma]{Corollary}

\newtheorem*{stat}{\name}
\newcommand{\name}{testing}
\theoremstyle{definition}

\newcommand{\bp}{\mathbin{\square}}

\newcommand{\bpz}{\boxdot}

\newcommand{\ltp}{\boxtimes}

\newcommand{\congtimes}{\mathbin{\square}}

\newcommand{\bboxtimes}{\mathbin{\ol{\congtimes}}}
\newcommand{\congtens}{\mathbin{\odot}}

\newcommand{\ootimes}{\mathbin{\ol{\otimes}}}
\newcommand{\jz}{$\set{\jj, 0}$}

\newcommand{\mL}[1]{M_3\langle#1\rangle}
\newcommand{\nL}[1]{N_5\langle#1\rangle}

\newenvironment{all}[1]{\renewcommand{\name}{#1}\begin{stat}}
                        {\end{stat}}

\DeclareMathOperator{\Con}{Con}
\DeclareMathOperator{\Aut}{Aut}
\DeclareMathOperator{\Pow}{Pow}
\DeclareMathOperator{\Conc}{Con_c}

\DeclareMathOperator{\Id}{Id}
\newcommand{\FL}{\mathrm{F}}
\newcommand{\T}{$(\mathrm{T})$}
\newcommand{\R}{\trianglelefteq}
\newcommand{\TJ}{$(\mathrm{T}_\jj)$}
\newcommand{\TM}{$(\mathrm{T}_\mm)$}

\newcommand{\ootimest}{\mathbin{\vec{\otimes}}}

\DeclareMathOperator{\Hom}{Hom}
\DeclareMathOperator{\J}{J}

\begin{document}

\title[A survey of tensor products]{A survey of tensor products\\ and
related constructions\\ in two lectures}

\author[G.~Gr\"atzer]{George Gr\"atzer}
\thanks{The research of the first author was
         supported by the NSERC of Canada.}
 \address{Department of Mathematics\\
	   University of Manitoba\\
	   Winnipeg MN, R3T 2N2\\
	   Canada}
 \email{gratzer@cc.umanitoba.ca}
 \urladdr{http://server.maths.umanitoba.ca/homepages/gratzer.html/}

 \author[F.~Wehrung]{Friedrich Wehrung}
 \address{C.N.R.S.\\
          Universit\'e de Caen, Campus II\\
          D\'epartement de Math\'ematiques\\
          B.P. 5186\\
          14032 Caen Cedex\\
          France}
 \email{wehrung@math.unicaen.fr}
 \urladdr{http://www.math.unicaen.fr/\~{}wehrung}


 \begin{abstract}
We survey tensor products of lattices with zero and related constructions
focused on two topics: amenable lattices and box products.
 \end{abstract}

\maketitle

\begin{center}
\tbf{PART I. FIRST LECTURE:\\
AMENABLE LATTICES}
\end{center}

 \begin{abstract}
Let $A$ be a finite lattice. Then $A$ is \emph{amenable} ($A \otimes B$ is
a lattice, for \emph{every} lattice $B$ with zero) if{}f $A$ (as a
join-semilattice) is \emph{sharply transferable} (whenever $A$ has an
embedding $\gf$ into $\Id L$, the ideal lattice of a lattice $L$, then $A$
has an embedding $\gy$ into $L$ satisfying $\gy(x) \in \gf(x)$ and
$\gy(x) \nin \gf(y)$, if $y < x$).

In Section \ref{S:Tensor}, we survey tensor products. In Section
\ref{S:Transferable}, we introduce transferability. These two topics are
brought together in Section~\ref{S:Amenable} in the characterization
theorem of amenable lattices.
 \end{abstract}

\section{Tensor product}\label{S:Tensor}
For a \jz-semilattice $A$, we use the notation $A^-=A -
\set{0}$.

Tensor products were introduced in J. Anderson and N. Kimura \cite{AK68}
and G.\,A. Fraser \cite{gF76}. Let $A$ and $B$ be \jz-semilattices. We
denote by $A \otimes B$ the \emph{tensor product} of $A$ and $B$, defined
as the free \jz-semilattice generated by the set $A^- \times B^-$ and subject
to the relations
\begin{equation*}
 \begin{aligned}
  \vv<a,b_0> \jj \vv<a,b_1> &= \vv<a, b_0 \jj b_1>,\quad
   &\text{for }a \in A^-,\ b_0,\,b_1  \in B^-;\\
  \vv<a_0, b> \jj \vv<a_1, b> &= \vv<a_0 \jj a_1, b>,\quad
   &\text{for }a_0,\,a_1 \in A^-,\ b  \in B^-.
 \end{aligned}
\end{equation*}

\subsection{The set representation}\label{S:set}
Let $A$ and $B$ be \jz-semilattices. We introduce a partial binary
operation, the \emph{lateral join}, on $A \times B$: let $\vv<a_0,b_0>$,
$\vv<a_1,b_1> \in A \times B$; the \emph{lateral join} $\vv<a_0,b_0> \jj
\vv<a_1,b_1>$ is defined if $a_0 = a_1$ or $b_0 = b_1$, in which
case, it is the join, $\vv<a_0\jj a_1,b_0\jj b_1>$; that is,
\begin{equation*}
 \begin{aligned}
  \vv<a,b_0> \jj \vv<a,b_1> &= \vv<a, b_0 \jj b_1>,\quad
   &\text{for }a \in A,\ b_0,\,b_1 \in B;\\
  \vv<a_0, b> \jj \vv<a_1, b> &= \vv<a_0 \jj a_1, b>,\quad
   &\text{for }a_0,\,a_1 \in A,\ b \in B.
 \end{aligned}
\end{equation*}

A nonempty subset $I$ of $A \times B$ is a \emph{bi-ideal} of $A \times
B$, if it is hereditary, it contains
 \[
   \bot_{A,B} = (A \times \set{0})  \uu  (\set{0} \times B),
 \]
and it is closed under lateral joins.

The \emph{extended tensor product} of $A$ and $B$, denoted by $A \ootimes
B$, is the lattice of all bi-ideals of $A \times B$. It is easy to see that
it is an algebraic lattice. For $a \in A$ and $b \in B$, we
define $a \otimes b \in A \ootimes B$ by
 \[
   a \otimes b =\bot_{A,B}  \uu  \setm{\vv<x, y>  \in A \times B}
  {\vv<x, y>  \leq \vv<a, b>}
 \]
and call $a \otimes b$ a \emph{pure tensor}.  A pure tensor is a
principal (that is, one-generated) bi-ideal.

Now we can state the representation:

\begin{theorem}\label{T:setrep}
The tensor product $A \otimes B$ can be represented as the
\jz-sub\-semi\-lat\-tice of compact elements of $A \ootimes B$.
\end{theorem}

Let $a_0 \leq a_1$ in $A$ and $b_0 \geq b_1$ in $B$.  Then
\[
      (a_0 \otimes b_0) \jj (a_1 \otimes b_1)  =
    (a_0 \otimes b_0)  \uu  (a_1 \otimes b_1).
\]
Such an element is called a \emph{mixed tensor}.

A bi-ideal $I$ is \emph{capped}, if it a \emph{finite union} of pure
tensors; pure tensors and mixed tensors are the simplest examples.  A
tensor product $A \otimes B$ is \emph{capped}, if (in the set
representation) all its elements are capped bi-ideals.  It is easy to see
that a capped tensor product is always a lattice. (It is an open problem
whether the converse holds; we do not think so.)

\subsection{Representation by homomorphisms}\label{S:homomorphisms}
Let $A$ and $B$ be \jz-semilattices. Note that $\Id B$, the set of all
ideals of $\vv<B; \jj>$, is a semilattice under intersection. So we can
consider the set of all semilattice homomorphisms from the semilattice
$\vv<A^-; \jj>$ into the semilattice $\vv<\Id B; \ii>$,
 \[
    A \ootimest B = \Hom(\vv<A^-; \jj>, \vv<\Id B; \ii>),
 \]
ordered componentwise, that is, $f\leq g$ if{}f
$f(a) \leq g(a)$ (that is,
$f(a) \ci g(a)$), for all $a \in A^-$. The arrow indicates which way the
homomorphisms go. Note that the elements of $A\ootimest B$ are
\emph{antitone} functions from $A^-$ to $\Id B$.

With any element $\gf$ of $A\ootimest B$, we associate the subset
$\ge(\gf)$ of $A \times B$:
 \[
  	\ge(\gf) = \setm{\vv<x, y> \in A \times B}{y \in \gf(x)} \uu
\bot_{A, B}.
 \]

\begin{theorem}\label{T:maprep}
The map $\ge$ is an isomorphism between $A\ootimest B$ and $A\ootimes
B$.
\end{theorem}

If $A$ is \emph{finite}, then a homomorphism from $\vv<A^-;\jj>$ to
$\vv<\Id B;\ii>$ is determined by its restriction to $\J(A)$, the set of
all join-irreducible elements of $A$.

For an interesting application of the representation of tensor products by
homomorphisms, see G. Gr\"atzer and F. Wehrung \cite{GW3}.

\subsection{Examples}\label{S:Examples}
 Let $B_n$ denote the Boolean lattice with $2^n$ elements.

Let $L$ be a lattice with zero.  Then
\begin{enumerate}
 \item $L \otimes B_1 \iso L$;
 \item $L \otimes B_n \iso L^n$;
 \item for a finite distributive lattice $D$ and $P = \J(D)$, $M_3
\otimes D$ can be represented as the set $M_3[D]$ of all balanced triples of
$D$ (a triple $\vv<x, y, z>$ is \emph{balanced} if{}f $x \mm y =x
\mm z =y \mm z$) or as $M_3^P$.
 \item $N_5 \otimes L$ can be represented as the set of all triples
$\vv<x, y, z>$ of $L$ satisfying $y \mm z \leq x \leq z$.
 \end{enumerate}

The representations in (iii) and (iv) utilize the representation by
homomorphisms of Section \ref{S:homomorphisms}.

The four examples share the property that the tensor product is a
lattice.  R.\,W. Quackenbush \cite{rQ85} raised the question whether this
is true, in general.  We answered this in \cite{GW0}.  In $M_3 \otimes
\textup{F}(3)$, let $a$, $b$, and $c$ be the atoms of $M_3$, let $x$, $y$,
and $z$ be the free generators of $\textup{F}(3)$, and form the elements
\begin{align*}
   \ga &= (a \otimes x) \jj (b \otimes y) \jj(c \otimes z),\\
   \gb &= a \otimes 1,
\end{align*}
where $1$ is the unit of $\FL(3)$. We proved that $\ga \mm \gb$ does not
exist in $M_3 \otimes \textup{F}(3)$.

\subsection{Congruences}\label{S:Congruences}
The main result of  G. Gr\"atzer, H. Lakser, and R.\,W. Quackenbush
\cite{GLQ81} is the statement that
 \[
   \Con A \otimes \Con B \iso \Con (A \otimes B)
 \]
holds for \emph{finite} lattices $A$ and $B$. For infinite lattices
with zero, this cannot hold, in general, because
\begin{itemize}
 \item the tensor product of two algebraic distributive lattices
is not necessarily algebraic;
 \item the tensor product of lattices with zero is not necessarily
a lattice.
\end{itemize}
 We compensate for the first by switching to the semilattice with zero of
compact congruences and for the second by assuming that the tensor
product is capped:

\begin{all}{The Isomorphism Theorem for Capped Tensor Products}
 Let $A$ and $B$ be lattices with zero. If $A \otimes B$ is capped, then
the following isomorphism holds:
 \[
   \Conc A \otimes \Conc B \iso \Conc (A \otimes B).
 \]
\end{all}

To describe this isomorphism, we need some notation. Let $\ga$ be a
congruence of $A$ and let $\gb$ be a congruence of $B$.  Define a binary
relation $\ga \bboxtimes \gb$ on $A  \ootimes  B$ as follows: for $H$, $K
\in A \ootimes  B$, let $\con H = K (\ga \bboxtimes \gb)$ if{}f,
for all $\vv<x, y>\in H$, there exists an $\vv<x', y'>\in K$ such
that $\con x = x'(\ga)$ and
$\con y = y'(\gb)$, and symmetrically. Let $\ga \congtimes \gb$ be the
restriction of $\ga \bboxtimes \gb$ to $A \otimes B$. If $A \otimes B$ is a
lattice, then $\ga \congtimes \gb$ is a lattice congruence on $A \otimes B$.

  For $\ga \in \Con A$ and $\gb \in \Con B$, we define  $\ga \congtens
\gb$, the \emph{tensor product} of $\ga$ and $\gb$, by the formula
 \[
   \ga \congtens \gb = (\ga \congtimes \go_B) \mm (\go_A \congtimes \gb).
 \]

\begin{theorem}\label{T:isomorphism}
Let $A$ and $B$ be lattices with zero such that $A \otimes B$ is a
lattice. The~map $\ga \otimes \gb \mapsto \ga \congtens \gb$ extends
to a \jz-embedding
 \[
   \ge \colon\Conc A \otimes \Conc B \to \Conc (A \otimes B).
 \]
If $A \otimes B$ is capped, then $\ge$ establishes the Isomorphism
Theorem.
\end{theorem}

The Isomorphism Theorem can be proved in a more general setup.

 Let $A$ and $B$ be lattices with zero. A \emph{sub-tensor product} of $A$
and $B$ is a subset $C$ of $A \otimes B$ satisfying the following
conditions:
 \begin{enumerate}
 \item $C$ contains all the mixed tensors in $A\otimes B$;
 \item $C$ is closed under finite intersection;
 \item $C$ is a lattice with respect to containment.
 \end{enumerate}

If every element of $C$ (as a bi-ideal) is capped, then $C$ is a
\emph{capped sub-tensor product}.

\begin{all}{The Isomorphism Theorem for Capped Sub-Tensor Products}
 Let $A$ and $B$ be lattices with zero. If $C$ is a capped sub-tensor
product of $A$ and $B$, then the following isomorphism holds:
 \[
   \Conc A \otimes \Conc B \iso \Conc C.
 \]
\end{all}

The lattice tensor product of Lecture Two is a sub-tensor product.

For some earlier results on congruence lattices of lattices of the type $L
\otimes D$, where $D$ is distributive, see  B. A. Davey, D. Duffus, R.\,W.
Quackenbush, and I.~Rival \cite{DQR}, D. Duffus, B. J\'onsson, and I.
Rival \cite{DJR}, J.\,D. Farley \cite{fF96},  G. Gr\"atzer and E.\,T.
Schmidt \cite{GS94}, G. Gr\"atzer and F. Wehrung \cite{GW4}, and E.\,T.
Schmidt \cite{tS79}.

\section{Transferable lattices}\label{S:Transferable}
Transferable lattices were introduced in \cite{gG70} in order to
provide a nice class of first-order sentences that hold for the ideal
lattice of a lattice if{}f they hold for the lattice.

 A finite lattice $T$ is \emph{transferable}, if for every embedding $\gf$
of $T$ into $\Id L$, the ideal lattice  of a lattice $L$, there exists an
embedding $\gx$ of $T$ into~$L$.

However, from a structural point of view, the following stronger form is
of more interest.

 A finite lattice $T$ is \emph{sharply transferable}, if for every
embedding $\gf$ of $T$ into $\Id L$, there exists an embedding $\gx$ of $T$
into~$L$ satisfying $\gx(x) \in \gf(y)$ if{}f $x \leq y$.

The motivation for these definitions comes from the fact that the
well-known result: \emph{a~lattice $L$ is modular if{}f\/ $\Id L$ is
modular}, can be recast: \emph{$N_5$ is a (sharply) transferable lattice}.

It is easy to verify that \emph{$N_5$ is a sharply transferable lattice}.
It is somewhat more difficult to see the negative result: \emph{$M_3$ is
not a (sharply) transferable lattice}.

To give the characterization theorem of (sharply) transferable lattices,
we need the following definitions, see H. Gaskill \cite{hG73}.

Let $P$ be a poset and let $X$ and $Y$ be subsets of $P$. Then $X$
\emph{is dominated by} $Y$, in notation, $X \ll Y$, if for all $x \in X$,
there exists $y \in Y$ such that $x \leq y$.

Let $A$ be a finite join-semilattice. A \emph{minimal pair} of $A$ is a
pair $\vv<p, I>$ such that $p \in \J(A)$, $I \ci \J(A)$, $|I| \geq 2$,
$p \nin I$, and $p \leq \JJ{I}$; moreover, for all $J \ci \J(A)$, if
$J \ll I$ and $p \leq \JJ{J}$, then $I \ci J$.

A finite join-semilattice $A$ satisfies condition \T, if $\J(A)$ has a linear
ordering $\R$ such that for every minimal pair $\vv<p, J>$ of $A$ and $j
\in J$, the relation $p \R j$ holds.  A~lattice $A$ satisfies condition
\TJ\ (respectively, \TM), if the semilattice $\vv<A; \jj>$ (respectively,
$\vv<A; \mm>$) satisfies \T.

Finally, we need the \emph{Whitman condition}:
 \begin{equation}
   \tag{W} x\mm y \leq u \jj v \text{\q implies that\q} [x \mm y, u
\jj v] \ii \set{x, y, u, v} \ne \es.
 \end{equation}

Now we can state the result from H.\,S. Gaskill, G. Gr\"atzer, and
C.\,R. Platt \cite{GGP75}:

\begin{all}{The Characterization Theorem for Sharply Transferable Lattices}
Let $A$ be a finite lattice. Then $A$ is sharply
transferable if{}f it satisfies the three conditions \TJ, \TM, and
\textup{(W)}.
\end{all}

As discussed in Appendix A and R. Freese's Appendix G of \cite{GLT2}, this
result shows that sharply transferable lattices are the same as finite
sublattices of a free lattice (see J.\,B. Nation \cite{jN82}).

Sharply transferable semilattices are defined analogously. H. Gaskill
\cite{hG73} proved the following result:

\begin{all}{The Characterization Theorem for Sharply Transferable
Semi\-lattices}\q \\
 Let $S$ be a finite semilattice. Then $S$ is sharply transferable if{}f it
satisfies \T.
\end{all}

See R.~Freese, J.~Je\v zek, and J.\,B. Nation \cite{FJN95} for a
discussion on how \TJ\ is the same as D-\emph{cycle free} and on the
structure of this class of lattices.

\section{Amenable lattices}\label{S:Amenable}
 Of course, the tensor product of two finite lattices is always a lattice.
In Section~\ref{S:Examples}, we noted that $M_3 \otimes \textup{F}(3)$ is
not a lattice.  Now we introduce the class of finite lattices $A$ for
which $A \otimes L$ is always a lattice.

Let us call the finite lattice $A$ \emph{amenable}, if $A \otimes L$ is
a lattice, for any lattice $L$ with zero.  So $M_3$ is not amenable. Every
finite distributive lattice is amenable. It is easy to see using the
representation in Example (iv) of Section~\ref{S:Examples} that $N_5$ is
amenable.

Now we state the characterization theorem of finite amenable lattices
\cite{GW1}:

\begin{theorem}\label{T:characterization_amenable}
 For a finite lattice $A$, the following
conditions are equivalent:
 \begin{enumerate}
\item $A$ is amenable.
\item $A$ is \emph{transferable} as a join-semilattice.
\item $A \otimes \FL(3)$ is a lattice.
\item $A$ satisfies \TJ.
\end{enumerate}
\end{theorem}

The equivalence of (i) and (iii) states that $\FL(3)$ is a ``test
lattice''; the
equivalence of (ii) and (iv) is a restatement of the result of H.~Gaskill
\cite{hG73} stated above.

The proof of this result is fairly long.  Curiously, the crucial step is
based on a construction in H.\,S. Gaskill, G. Gr\"atzer, and C.\,R. Platt
\cite{GGP75} for lattice (not semilattice) transferability; while we are
unable to apply this result directly, the idea is clearly borrowed.

It follows that the class of finite amenable lattices and the class of
finite \emph{lower bounded} lattices coincide, see R.~Freese, J.~Je\v zek,
and J.\,B. Nation \cite{FJN95}. By Theorem~2.43 of \cite{FJN95}, a
finite lattice is lower bounded if{}f it can be obtained from a
one-element lattice by a sequence of doubling constructions with
respect to lower pseudo-intervals.

Recently, we have succeeded in generalizing
Theorem~\ref{T:characterization_amenable} to arbitrary lattices
with zero:
 
 \begin{theorem}\label{T:characterization_amenable_general}
For a lattice $A$ with zero, the following conditions are
equivalent:
 \begin{enumerate}
\item $A$ is amenable.
\item $A$ is locally finite and
$A\otimes B$ is a lattice, for every lattice $B$ with zero.
\item $A$ is locally finite and $A \otimes \FL(3)$ is a lattice.
\item $A$ is locally finite and every finite sublattice of
$A$ satisfies \TJ.
\end{enumerate}
\end{theorem}

For a finite amenable lattice $A$, there is a close connection between
$\J(A)$ and $\J(\Con A)$.  Let $a \in \J(A)$; let $a_*$ be the unique
element of $A$ covered by $a$.  Then $a \mapsto \gQ(a, a_*)$ is a
bijection between $\J(A)$ and $\J(\Con A)$.  (In fact, the converse is
also true, showing that amenability is the same as \emph{fermentability} in the
sense of
 P. Pudl\'ak and J. T\r{u}ma \cite{PT74}.)  This suggests that the congruence
lattice of a finite amenable lattice is very special.

A \emph{spike} in a finite poset $P$ is a pair $a < b$ of elements of $P$
such that $b$ is maximal in $P$, $b$ covers $a$ in $P$, and $b$ is the
\emph{only} maximal element of $P$ containing $a$.  A poset $P$ is
\emph{spike-free}, if it has no spikes.

\begin{theorem}\label{T:spike-free}
 A finite distributive lattice $D$ can be represented as the congruence
lattice of an amenable lattice if{f} ${\textup J}(D)$ is spike-free.
\end{theorem}

This result is a special case of a more general theorem in \cite{GW5}.

\newpage

\begin{center}
\tbf{PART II. SECOND LECTURE:\\
BOX PRODUCTS}
\end{center}

 \begin{abstract}
 We have seen in Part I that the tensor product of two lattices with zero is
not necessarily a lattice. We survey a new lattice construction, the
\emph{box product} that always yields a lattice. If $A$ and $B$ are
lattices and either both $A$ and $B$ have a zero or one of them is bounded,
then the box product $A\bp B$ of $A$ and $B$ has an ideal, $A\ltp B$, for
which an analogue of the Isomorphism Theorem for capped sub-tensor products
holds, without any further restriction on $A$ or~$B$. In general, $A\ltp B$
is a subset of $A\otimes B$; equality holds, if $A$ or~$B$ is distributive.
 \end{abstract}

\section{The $\mL{L}$ construction and the $\nL{L}$ construction}
\label{S:M3N5L}

Let $L$ be a lattice. A lattice $K$ is a \emph{congruence-preserving
extension} of $L$, if $K$ is an extension of $L$ and every congruence of
$L$ extends to exactly one congruence of $K$. The extension is
\emph{proper}, if $K\ne L$. Similarly, we can define a congruence-preserving
embedding of lattices. In \cite{GS95}, the first author and E.\,T. Schmidt asked
whether every lattice $L$ with more than one element has a proper
congruence-preserving extension. If $L$ is a modular lattice, the answer is
already provided by Schmidt's $M_3[L]$ construction, see E.\,T. Schmidt
\cite{tS68}, R.\,W. Quackenbush \cite{rQ85}, and
Section~\ref{S:Examples}. By definition, $M_3[L]$ is the set of all
\emph{balanced triples} of $L$, ordered componentwise, see
Section~\ref{S:Examples}:
 \[
 M_3[L]=\setm{\vv<x,y,z>\in L^3}{x\mm y=x\mm z=y\mm z}.
 \]

Unfortunately, $M_3[L]$ is not always a lattice, see G. Gr\"atzer and F.
Wehrung
\cite{GW4} for a planar example $L$. The answer to the problem mentioned in the
previous paragraph was finally provided by a simple trick that we describe now,
see \cite{GW}. For every lattice $L$, define $\mL{L}$, a subset of $L^3$, as
follows:
 \begin{equation}\label{Eq:Bool}
   \mL{L} = \setm{\vv<v \mm w, u \mm w, u \mm v>}{u,\, v,\, w \in L}.
 \end{equation}
 We call an element of $\mL{L}$ a \emph{Boolean triple} of $L$. In particular,
$\mL{L}$ is a subset of $M_3[L]$. Endow $\mL{L}$ with the componentwise
ordering.

\begin{theorem}\label{T:M3new}
Let $L$ be a lattice. Then $\mL{L}$ is a lattice, and the diagonal map,
 \[
    x \mapsto \vv<x, x, x>,
 \]
defines a congruence-preserving embedding from $L$ into $\mL{L}$.
\end{theorem}

In particular, if $L$ has more than one element, then $\mL{L}$
properly contains $L$, thus solving the above problem.

It appears desirable to generalize the $\mL{L}$ construction to any pair of
lattices with zero, thus creating an analogue of the tensor product that never
fails to be a lattice. One (heuristic) way to proceed is the following. We note
that the Boolean triples of $L$ are exactly those triples of $L$ that are
balanced ``for a good reason''. Of course, one has to define precisely what a
``good reason'' is. Formula \eqref{Eq:Bool} suggests to look for
``meet-parametrizations''
of the solutions of the equational system defining
balanced triples, that is, $x \mm y = x \mm z = y \mm z$.

Now let us do the same with the pentagon, $N_5$, instead of $M_3$. By
using the representation by homomorphisms of the elements of the tensor
product $N_5\otimes L$, see Section~\ref{S:homomorphisms}, we define a
certain object that we denote by $N_5[L]$, see Section~\ref{S:Examples}:
 \begin{equation}\label{Eq:N5L}
 N_5[L]=\setm{\vv<x,y,z>\in L^3}{y\mm z\leq x\leq z}.
 \end{equation}
The situation here is quite different from the situation with $M_3[L]$:
indeed, since $N_5$ is amenable, $N_5[L]$ is always a lattice; furthermore,
if $L$ has a zero, then $N_5[L]$ is isomorphic to $N_5\otimes L$.

However, we may still look
for those triples of elements of $L$ that
belong to $N_5[L]$ ``for a good reason'' (say, a meet-parametrization of
the solutions of the equational system defining $N_5[L]$). An easy
computation gives us the definition of a new object that we denote, of
course, by $\nL{L}$:
 \begin{equation}\label{Eq:N5newL}
 \nL{L}=\setm{\vv<v\mm w,u\mm w,v>}{u,\ v,\ w\in L}.
 \end{equation}
Again, it is not hard to prove that $\nL{L}$, endowed with componentwise
ordering, is a lattice. It is strange that even though $N_5[L]$ is a
lattice, for every lattice $L$, $\nL{L}$ is, as a rule, a \emph{proper}
subset of $N_5[L]$; for example, for $L = N_5$.

A similar method to the one outlined above gives a definition of
$A\langle L\rangle$, for a finite lattice $A$ and a lattice $L$. A precise
description of this method would be lengthy, and it would involve
the study of the structure of solution sets of systems of equations in
distributive semilattices. Furthermore, it may not be very useful at this
point, because we found a general, short definition that encompasses all
these constructions and more. The starting point is the construction of the
\emph{box product} defined in the next section.

\section{The box product $A\bp B$}\label{S:BoxProd}

We refer to \cite{GW2}, for more detail and for proofs.

 Let $A$ and $B$ be lattices. For $\vv<a,b>\in A\times B$, define
 \[
    a \bp b = \setm{\vv<x, y> \in A \times B}{x \leq a \text{ or }y \leq b}.
 \]
We define the \emph{box product} of $A$ and $B$, denoted by $A \bp B$, as
the set of all finite \emph{intersections} of the form
 \[
    H = \IIm{a_i \bp b_i}{i < n},
 \]
where $n$ is a positive integer and $\vv<a_i, b_i> \in A \times B$, for all
$i < n$.

It is clear that $A\bp B$ is a meet-subsemilattice of
the powerset lattice $\Pow(A\times B)$ of $A\times B$. To obtain that
$A\bp B$ is also a join-semilattice, we prove that it is a closure system
in a sublattice, denoted by $A\bpz B$, of $\Pow(A\times B)$. The
definition of $A\bpz B$ is the following. For $\vv<c,d>\in A\times B$, put
 \[
 c\circ d=\setm{\vv<x,y>\in A\times B}{x\leq c\text{ and }y\leq d},
 \]
and define $A\bpz B$ as the set of all finite \emph{unions} of the form
 \begin{equation}\label{Eq:Hbpz}
 H=\UUm{a_i\bp b_i}{i<m}\uu\UUm{c_j\circ d_j}{j<n},
 \end{equation}
where $m>0$, $n\geq0$, and all pairs $\vv<a_i,b_i>$ and $\vv<c_j,d_j>$
belong to $A\times B$.

\begin{theorem}\label{T:BoxProd}
Let $A$ and $B$ lattices. Then $A\bpz B$ is a sublattice of
$\Pow(A\times B)$ and $A\bp B$ is a closure system in $A\bpz B$. In
particular, $A\bp B$ is a lattice.
\end{theorem}

The statement that $A \bp B$ is a closure system in $A \bpz B$ means that,
for every element $H$ of $A\bpz B$, there exists a least element $K$
of $A\bp B$ such that $H\ci K$; we denote this element by $\ol{H}$. It is
important to note that $\ol{H}$ is given by a \emph{formula}, as follows.
If $H$ is written as in \eqref{Eq:Hbpz}, then $\ol{H}$ is given by
 \[
   \ol{H}=\IIm{a^{(X)}\bp b^{(n-X)}}{X\ci n},
 \]
where
 \begin{align*}
   a^{(X)}=\JJm{a_i}{i<m}\jj\JJm{c_j}{j\in X},\\
   b^{(X)}=\JJm{b_i}{i<m}\jj\JJm{d_j}{j\in X},
 \end{align*}
for all $X\ci n$.

\section{The lattice tensor product $A\ltp B$}

For lattices $A$ and $B$, the box product $A\bp B$ has a unit
element if{}f
either $A$ or $B$ has a unit element. In particular, $M_3\bp B$ always has a
unit element, so that it is not isomorphic to the lattice $\mL{B}$ of
Boolean triples of $B$, see Section~\ref{S:M3N5L}. Thus we shall define an
\emph{ideal} of $A\bp B$.

For arbitrary lattices $A$ and $B$, we can modify the definition of
$\bot_{A,B}$, introduced in Section~\ref{S:Tensor}, as follows:
 \[
    \bot_{A,B}=(A\times\bot_B)\cup(\bot_A\times B),
 \]
 where
 \begin{equation*}
   \bot_L =
           \begin{cases}
              \set{0_L}, &\text{if $L$ has a zero,}\\
              \es,       &\text{otherwise.}
          \end{cases}
 \end{equation*}
For $\vv<a,b>\in A\times B$, define
 \[
 a\ltp b=\setm{\vv<x,y>\in A\times B}{x\leq a\text{ and }y\leq b}
 \cup\bot_{A,B}.
 \]
If both $A$ and $B$ have a zero element, then $a\ltp b$ is an element of
$A\bp B$, namely,
 \[
 a\ltp b=(a\bp 0_B)\ii(0_A\bp b)=a\otimes b.
 \]

 An element $H$ of $A \bp B$ is \emph{confined}, if it is contained in
$a\ltp b$ for some $\vv<a,b>\in A\times B$. We define $A\ltp B$, the
\emph{lattice tensor product} of $A$ and $B$, as the ideal of $A\bp B$
consisting of all confined elements of $A\bp B$.

If $A$ has a zero element and $B$ has no zero element, then $a\ltp b$ does
not contain any element of $A\bp B$ unless $A$ has a unit element (so that
$A$ is bounded), in which case $a$ equals this unit, thus $a\ltp b$ equals
$0_A\bp b$, so that it belongs to $A\bp B$. In particular, if $A$ has a zero
but no unit and $B$ has no zero, then $A\ltp B=\es$. In fact, it is easy to
see exactly when $A\ltp B$ is nonempty:

\begin{lemma}\label{L:AltpBnes}
Let $A$ and $B$ be lattices. Then $A\ltp B$ is nonempty if{}f
one of the
following conditions holds:
 \begin{enumerate}
 \item both $A$ and $B$ have zero;
 \item either $A$ or $B$ is bounded;
 \item both $A$ and $B$ have unit.
 \end{enumerate}
 \end{lemma}

In case (iii), that is, if both $A$ and $B$ have unit, then every element of
$A\bp B$ is bounded, so that $A\ltp B=A\bp B$. For a lattice $L$, denote by
$L^{\mathrm{d}}$ the dual lattice of~$L$.

As one would expect, cases (i) and (iii) correspond to each other \emph{via}
lattice dualization:

\begin{theorem}\label{T:bpdualltp}
Let $A$ and $B$ be lattices with zero. Then the following isomorphism holds:
 \[
   (A\ltp B)^{\mathrm{d}}\iso A^{\mathrm{d}}\bp B^{\mathrm{d}}.
 \]
 \end{theorem}

Interestingly, the main observation on the Isomorphism Theorem for
lattice tensor products concerns lattices \emph{with zero} (as opposed to
lattices with unit):

 \begin{theorem}\label{T:ltpcapped}
 Let $A$ and $B$ be lattices with zero. Then $A \ltp B$ is a capped sub-tensor
product of $A$ and $B$. Furthermore, $A\ltp B$ is the \emph{smallest} capped
sub-tensor product of $A$ and $B$, with respect to containment.
 \end{theorem}

The Isomorphism Theorem for Capped Sub-Tensor Products, see
Section~\ref{S:Congruences}, implies then that the isomorphism $\Conc(A \ltp B)
\iso \Conc A \otimes \Conc B$ holds, for lattices $A$ and $B$ with zero. A
direct
limit argument and some extra work makes it then possible to obtain the
following
general result:

\begin{theorem}\label{T:Isoltp}
Let $A$ and $B$ be lattices. If $A\ltp B$ is nonempty, then the
following isomorphism holds:
 \[
 \Conc(A\ltp B)\iso\Conc A\otimes\Conc B.
 \]
\end{theorem}

Theorem~\ref{T:Isoltp} is proved by constructing a map,
 \[
   \gm\colon\Conc A\otimes\Conc B\longrightarrow\Conc(A\ltp B),
 \]
 and proving that $\gm$ is an isomorphism. The isomorphism $\gm$ is easy to
describe. Since $\gm$ is a join homomorphism, it is sufficient to describe the
image of a pure tensor $\ga\otimes\gb$, where $\ga=\gQ_A(a_0, a_1)$ and
$\gb=\gQ_B(b_0, b_1)$ (with $a_0 \leq a_1$ in $A$ and $b_0 \leq b_1$ in $B$).
According to Lemma~\ref{L:AltpBnes}, we split the description into three cases:

 \begin{enumerate}
 \item $A$ and $B$ are lattices with zero:
 \[
   \gm(\ga\otimes\gb)=
    \gQ_{A\ltp B}((a_0\ltp b_1)\jj(a_1\ltp b_0),a_1\ltp b_1).
 \]
 \item $A$ is bounded (or symmetrically, $B$ is bounded):
 \[
   \gm(\ga\otimes\gb)=
     \gQ_{A\ltp B}((a_0\bp b_0)\ii(0_A\bp b_1),(a_1\bp b_0)\ii(0_A\bp b_1)).
 \]
 \item $A$ and $B$ are lattices with unit:
 \[
   \gm(\ga\otimes\gb)=
    \gQ_{A\ltp B}(a_0\bp b_0,(a_0\bp b_1)\ii(a_1\bp b_0)).
 \]
\end{enumerate}

Of course, formula (iii) can be obtained from formula (i) and the
canonical isomorphism given in Theorem~\ref{T:bpdualltp}.

The lattice tensor product construction $A\ltp B$ can be easily related to
the constructions $\mL{L}$ and $\nL{L}$ described in Section~\ref{S:M3N5L}:

\begin{theorem}\label{T:M3N5ltp}
Let $L$ be a lattice. Then the following isomorphisms hold:
 \begin{align*}
   M_3 \ltp L &\iso \mL{L},\\
   N_5 \ltp L &\iso \nL{L}.
 \end{align*}
 \end{theorem}

An isomorphism $\ga \colon \mL{L} \to M_3 \ltp L$ is given by
 \[
   \ga(\vv<v\mm w,u\mm w,u\mm v>)=(p\bp u)\ii(q\bp v)\ii(r\bp w),
 \]
for all $u$, $v$, $w\in L$, where $p$, $q$, and $r$ are the atoms of
$M_3$.

 An isomorphism $\gb \colon \nL{L} \to N_5 \ltp L$ is given by
 \[
   \gb(\vv<v\mm w,u\mm w,v>)=(a\bp u)\ii(b\bp v)\ii(c\bp w),
 \]
where $a>c$ and $b$ are the join-irreducible elements of $N_5$.
\smallskip

Much more general is the following corollary of Theorem~\ref{T:Isoltp} and
of the formulas describing the isomorphism $\gm$:\goodbreak

\begin{corollary}\label{C:Isoltp}
Let $S$ and $L$ be a lattices; let $S$ be simple.
 \begin{enumerate}
 \item If $S$ is bounded, then the map $j\colon L\to S\ltp L$ defined by
 \[
   j(x)=0_S\bp x,
 \]
 for all $x\in L$, is a congruence-preserving lattice embedding.

\item If both $S$ and $L$ have zero, then for every $s\in S^-$, the map
$j_s\colon L\to S\ltp L$ defined by
 \[
   j_s(x)=s\ltp x,
 \]
for all $x\in L$, is a congruence-preserving lattice embedding.
\end{enumerate}
\end{corollary}

For $S=M_3$ and \emph{via} the identification of $M_3\ltp L$ with
$\mL{L}$, the first embedding is the map $x\mapsto\vv<x,x,x>$, while the
second embedding is, for example, for $s=p$, the map $x\mapsto\vv<x,0,0>$.
For more general $S$, this can be used to prove statements stronger than
Theorem~\ref{T:M3new}, such as the \emph{Strong Independence Theorem}, see
Section~\ref{S:StrInd}.

\section{Some applications}\label{S:LE}

\subsection{Congruence representations of distributive semilattices with
zero}\label{S:Repr}
 Let us say that a \jz-semilattice $S$ is \emph{representable}
(\emph{$\set{0}$-rep\-re\-sent\-a\-ble},
\emph{$\set{0,1}$-rep\-re\-sent\-a\-ble}, respectively), if there exists a
lattice
$L$ (a lattice $L$ with zero, a bounded lattice $L$, respectively) such
that $\Conc L\iso
S$. It is an open problem, dating back to the forties, whether every
distributive
\jz-semilattice is representable or $\set{0}$-rep\-re\-sent\-a\-ble. Similarly,
it is an open problem whether every bounded distributive \jz-semilattice is
representable, or $\set{0}$-rep\-re\-sent\-a\-ble, or
$\set{0,1}$-rep\-re\-sent\-a\-ble. We refer to G. Gr\"atzer and E.\,T. Schmidt
\cite{GS98} for a detailed history of this problem. We recall here some partial
answers:

\begin{enumerate}
\item If $S$ satisfies one of the following conditions, then $S$ is
representable
(see \cite{GS98}, Theorem 13): \begin{enumerate}
\item $\Id S$ is completely distributive (R.~P. Dilworth);

\item $S$ is a lattice (E.\,T. Schmidt);

\item $S$ is locally countable, that is, every element of $S$ generates a
countable principal ideal (A.~P. Huhn for $S$ countable, H. Dobbertin in
general).

\item $|S|\leq\aleph_1$ (A.~P. Huhn).
\end{enumerate}

In all four cases, the representability of $S$ can be obtained \emph{via} E.\,T.
Schmidt's condition (see \cite{tS68}) that $S$ is a distributive image of a
generalized Boolean semilattice. A closer look at the proofs shows that, in
fact,
Schmidt's condition implies $\set{0}$-representability.

 \item If $S$ is countable, then $S$ is representable by a \emph{sectionally
complemented modular lattice} $L$ (G.\,M. Bergman \cite{Berg86}, see also
K.\,R. Goodearl and F. Wehrung \cite{GoWe}). Furthermore, if $S$ is bounded,
then one can take $L$ to be bounded.

 \item If $|S|\leq\aleph_1$, then $S$ is representable by a relatively
complemented (not modular \emph{a priori}) lattice with zero. The proof of this
result is based on an amalgamation result of J. T\r uma \cite{Tuma}, see
also G.
Gr\"atzer, H. Lakser, and F.~Wehrung \cite{GLW}. However, the method fails to
produce a bounded lattice $L$ even if $S$ is bounded.

\end{enumerate}

New consequences can be obtained about the class $\E R$ of representable
\jz-semilattices, the class $\E R_0$ of $\set{0}$-rep\-re\-sent\-a\-ble
\jz-semilattices and the class $\E R_{0,1}$ of
$\set{0,1}$-rep\-re\-sent\-a\-ble \jz-semilattices, by using
Theorem~\ref{T:Isoltp}:\goodbreak

\begin{corollary}\label{C:RepClo}\hfill
\begin{enumerate}
\item The classes $\E R_0$ and $\E R_{0,1}$ are closed under tensor
product.

\item Let $A\in\E R_{0,1}$ and let $B\in\E R$. Then $A\otimes B\in\E R$.
\end{enumerate}
\end{corollary}

This result can be extended to \emph{iterated tensor products}. If $\vv<S_i\mid
i\in I>$ is a family of bounded \jz-semilattices, then their \emph{iterated
tensor product} is the direct limit of the family $\bigotimes_{i\in J}S_i$,
where
$J$ ranges over all finite subsets of $I$, and the transition homomorphisms are
defined by $\otimes_{i\in J}x_i\mapsto\otimes_{i\in K}x_i$, where
$x_i=1_{S_i}$,
for $i\in K-J$, and $J \ci K$ are finite subsets of $I$.

 \begin{corollary}\label{C:RepCloIt}
 The class $\E R_{0,1}$ is closed under iterated tensor products.
 \end{corollary}

Further results can be obtained for other subclasses of $\E R$. Let us
mention, for example, the following. If $L$ is a lattice, we say that $L$
has \emph{permutable congruences}, if any two congruences of $L$
commute.

\begin{lemma}\label{L:CommCong}
Let $A$ and $B$ be lattices such that $A\ltp B$ is nonempty. If $A$
and $B$ have permutable congruences, then $A\ltp B$ has permutable
congruences.
\end{lemma}

By the known representation results, the class of all \jz-semilattices that are
representable by lattices with zero and with permutable congruences
contains all distributive semilattices of size at most
$\aleph_1$---this is because every relatively complemented lattice has
permutable congruences. Denote by $\E R^{\mathrm{c}}$ ($\E
R_0^{\mathrm{c}}$, $\E R_{0,1}^{\mathrm{c}}$, respectively) the
class of all \jz-semilattices that are representable by lattices (lattices with
zero, bounded lattices, respectively) with permutable congruences. It
is proved in J. T\r uma
and F. Wehrung \cite{TuWe}, using the main result of  M. Plo\v s\v cica, J.
T\r uma and F. Wehrung \cite{PTWe}, that $\E R^{\mathrm{c}}$ is a proper
subclass of $\E R$.

\begin{corollary}\label{C:RepCloComm}\hfill
\begin{enumerate}
\item The classes $\E R_0^{\mathrm{c}}$ and $\E R_{0,1}^{\mathrm{c}}$ are
closed under tensor product.

\item Let $A\in\E R_{0,1}^{\mathrm{c}}$ and let $B\in\E R^{\mathrm{c}}$.
Then $A\otimes B\in\E R^{\mathrm{c}}$.

\item The class $\E R_{0,1}^{\mathrm{c}}$ is closed under iterated tensor
product.
\end{enumerate}
\end{corollary}

There is an intriguing similarity between these preservation results and
known representation results of dimension groups as ordered $K_0$ groups of
locally matricial rings, see K.\,R. Goodearl and D.\,E. Handelman
\cite{GoHa}.

\subsection{Strong independence of the congruence lattice and the
automorphism group}\label{S:StrInd}
The \emph{Independence Theorem} for the congruence
lattice and the automorphism group of a finite lattice
was proved by V.\,A. Baranski\u\i\ \cite{vB79} and A.~Urquhart
\cite{aU78} (solving Problem II.19 of \cite{GLT1}) :

\begin{all}{The Independence Theorem for Finite Lattices}
Let $G$ be a finite group and let $D$ be a finite distributive lattice. Then
there exists a finite lattice $L$ such that $\Aut L$, the automorphism group of
$L$, is isomorphic to $G$, while $\Con L$, the congruence lattice of $L$, is
isomorphic to $D$.
 \end{all}

Both proofs utilize the characterization theorem of
congruence lattices of finite lattices (as finite
distributive lattices) and the characterization
theorem of automorphism groups of finite lattices (as
finite groups).

In G. Gr\"atzer and E.\,T. Schmidt \cite{GS95b}, a new,
stronger form of independence is introduced.

A finite lattice  $K$ is an \emph{automorphism-preserving extension} of
$L$, if
$K$ is an extension  and every automorphism of $L$ has exactly one extension to
$K$, and in addition, every automorphism of $K$ is the extension of an
automorphism of $L$.  Of course, then the automorphism group of $L$ is
isomorphic
to the automorphism group of $K$.

The following result has been established in G. Gr\"atzer and E.\,T. Schmidt
\cite{GS95b}:

\begin{all}{The Strong Independence Theorem for Finite Lattices}
Let $L_{\mathrm{C}}$ and $L_{\mathrm{A}}$ be finite lattices,
let $L_{\mathrm{C}}$ have more than one element, and let
$L_{\mathrm{C}} \ii L_{\mathrm{A}} = \set{0}$. Then~there exists a finite
atomistic lattice $L$ that is a congruence-preserving extension of
$L_{\mathrm{C}}$ and an automorphism-preserving extension
of~$L_{\mathrm{A}}$. In fact, both extensions preserve the zero.
\end{all}

Of course, the congruence lattice of $L$ is isomorphic to the congruence
lattice
of $L_{\mathrm{C}}$, and the automorphism group of $L$ is isomorphic to the
automorphism group of  $L_{\mathrm{A}}$. Therefore, indeed, for finite
lattices,
independence follows from strong independence. This is because every finite
distributive lattice can be obtained as $\Con L_{\mathrm{C}}$ for some finite
lattice $L_{\mathrm{C}}$ (R.~P. Dilworth; see G. Gr\"atzer and E.\,T. Schmidt
\cite{GS62}) and every finite group can be obtained as
$\Aut L_{\mathrm{A}}$ for
some finite lattice $L_{\mathrm{A}}$ (see G. Birkhoff \cite{Bi}).

The question of a possible generalization of the Independence Theorem or the
Strong Independence Theorem to infinite lattices was raised in Problems~1
and~2 of G. Gr\"atzer and E.\,T. Schmidt \cite{GS95b} (Problem 3, whether
every lattice with more than one element has a proper congruence-preserving
extension, is solved in our paper \cite{GW}, see Theorem~\ref{T:M3new}). The
statement of independence for arbitrary lattices is by itself a problem,
because
it is not known which distributive \jz-semilattices $S$ are representable as
$\Con L$ for a lattice $L$---that is, which $S$ belong to the class $\E R$, see
Section~\ref{S:Repr}. On the other hand, Birkhoff's result extends to all
groups:
\emph{every group is isomorphic to the automorphism group of some lattice}.
Thus
a possible formulation of independence for infinite lattices would be with
\emph{representable} \jz-semilattices, on the one hand, and arbitrary
groups, on
the other. Again, such a statement would follow from strong independence.

We proved strong independence in G. Gr\"atzer and F. Wehrung \cite{GW6}, thus
solving Problem~II.18 of \cite{GLT1} and Problems 1 and 2 of \cite{GS95b}:

 \begin{all}{The Strong Independence Theorem for Lattices with Zero}
 \ Let $L_{\mathrm{A}}$ and $L_{\mathrm{C}}$ be lattices with zero, let
$L_{\mathrm{C}}$ have more than one element.
Then there exists a lattice $L$ that is a $\set{0}$-preserving extension of
both $L_{\mathrm{A}}$ and $L_{\mathrm{C}}$, an automorphism-preserving
extension of $L_{\mathrm{A}}$, and a congruence-preserving extension of
$L_{\mathrm{C}}$.
 \end{all}

 \begin{all}{The Strong Independence Theorem for Lattices}
 Let $L_{\mathrm{A}}$ and $L_{\mathrm{C}}$ be lattices, let
$L_{\mathrm{C}}$ have more than one element.
Then there exists a lattice $L$ that is an automorphism-preserving extension of
$L_{\mathrm{A}}$ and a congruence-preserving extension of $L_{\mathrm{C}}$.
 \end{all}

The main ingredients of the proof are direct limits, gluings, and box
products (in fact, lattice tensor products).

\end{document}